\documentclass[12pt,amsfonts]{article}
\usepackage{graphicx}
\usepackage{latexsym}
\usepackage{amssymb}
\usepackage{amsmath}
\usepackage{enumerate}

\usepackage{layout}
\usepackage{eufrak}
\newtheorem{prop}{Proposition}
\newtheorem{lemma}{Lemma}

\newtheorem{corollary}{Corollary}
\newtheorem{theorem}{Theorem}
\newtheorem{remark}{Remark}
\newtheorem{example}{Example}

\def\real{{\mathord{{\rm I\kern-2.8pt R}}}}        
\def\inte{{\mathord{{\rm I\kern-2.8pt N}}}}

\def\sZZ{{\rm Z\kern-2.8ptem{}Z}}

\def\z{{\mathchoice
  {\sZZ}
  {\sZZ}
  {\rm Z\kern-0.30em{}Z}
  {\rm Z\kern-0.25em{}Z} }}
\def\sQQ{{\kern 0.27em \vrule height1.45ex width0.03em depth0em
          \kern-0.30em \rm Q}}
\def\qu{{\mathchoice
    {\sQQ}
    {\sQQ}
  {\kern 0.225em \vrule height1.05ex width0.025em depth0em \kern-0.25em \rm Q}
  {\kern 0.180em \vrule height0.78ex width0.020em depth0em \kern-0.20em \rm Q}
        }}
\def\sCC{{\kern 0.27em \vrule height1.45ex width0.03em depth0em
          \kern-0.30em \rm C}}
\def\complex{{\mathchoice
    {\sCC}
    {\sCC}
  {\kern 0.225em \vrule height1.05ex width0.025em depth0em \kern-0.25em \rm C}
  {\kern 0.180em \vrule height0.78ex width0.020em depth0em \kern-0.20em \rm C}
        }}

\newcommand{\ba}{\begin{array}}
\newcommand{\ea}{\end{array}}
\newcommand{\be}{\begin{equation}}
\newcommand{\ee}{\end{equation}}
\newcommand{\bea}{\begin{eqnarray}}
\newcommand{\eea}{\end{eqnarray}}
\newcommand{\beaa}{\begin{eqnarray*}}
\newcommand{\eeaa}{\end{eqnarray*}}

\def\z{\zeta}

\font\tenmath=msbm10 \font\sevenmath=msbm7 \font\fivemath=msbm5
\newfam\mathfam \textfont\mathfam=\tenmath
\scriptfont\mathfam=\sevenmath \scriptscriptfont\mathfam=\fivemath

\def \={{\buildrel {\rm (law)} \over =}}

\def\qed{ \hfill \vrule width.25cm height.25cm depth0cm\smallskip}

\newcommand{\basa}{\begin{assumption}}
\newcommand{\easa}{\end{assumption}}

\newcommand{\bas}{\begin{assum}}
\newcommand{\eas}{\end{assum}}

\newcommand{\ignore}[1]{}
\begin{document}

\renewcommand{\thefootnote}{\fnsymbol{footnote}}

\renewcommand{\thefootnote}{\fnsymbol{footnote}}

\title{The determinant of the Malliavin matrix and the determinant of the covariance matrix for multiple integrals}

\author{Ciprian A. Tudor \footnote{  Supported by the CNCS grant PN-II-ID-PCCE-2011-2-0015. }\\
 Laboratoire Paul Painlev\'e, Universit\'e de Lille 1\\
 F-59655 Villeneuve d'Ascq, France\\
 and\\
Academy of Economical Studies, Bucharest, Romania\\
\quad tudor@math.univ-lille1.fr}
\maketitle
 \maketitle

\begin{abstract}
A well-known problem in Malliavin calculus concerns the relation between the determinant of the Malliavin matrix of a random vector and the determinant of its covariance matrix. We give an explicit relation between these two determinants for couples of random vectors of multiple integrals. In particular, if the multiple integrals are of the same order and this order is at most 4, we prove that two random variables in the same Wiener chaos either admit a joint density, either are proportional and that the result is not true for random variables in Wiener chaoses of different orders. 
\end{abstract}

\vskip0.3cm

{\bf 2010 AMS Classification Numbers:}   60F05, 60H05, 91G70.

 \vskip0.3cm


\section{Introduction}
The original motivation of the Malliavin calculus was to study the existence and the regularity of the densities of random variables.  In this research direction, the determinant of the so-callled Malliavin matrix plays a crucial role. 

We give here an explicit formula that connects the determinant of the Malliavin matrix and the determinant of the covariance matrix of a couple of multiple stochastic integrals. This is related to two open problems stated in \cite{NoNuPo}. In this reference, the authors showed that, if $F=(F_{1},..,F_{d}) $ is a random vector whose components belong to a finite sum of Wiener chaoses, then the law of $F$ is not absolutely  continuous with respect to the Lebesque measure if and only if $E \det \Lambda =0$. Here $\Lambda $ denotes the Malliavin matrix of the vector $F$. In particular, they proved that  a couple of multiple integrals of order 2 either admits a density or its components are proportional.

They stated two open questions (Questions 6.1 and 6.2 in \cite{NoNuPo}, arXiv version): if $C$ is the covariance matrix and $\Lambda$ the Malliavin matrix of a vector of multiple stochastic integrals, 
\begin{description}
\item{$\bullet$ } is there true that $E\det \Lambda \geq c \det C$, with $c>0$ an universal constant?

\item{$\bullet$ } is there true that the law of  a vector of multiple integrals with components in the same Wiener chaos is either absolutely continuous with respect to the Lesque measure or its components are proportional?
\end{description}

We make a first step in order to answer to  these two open problems. Actually, we find an explicit relation that connects the two determinants. In particular, if the multiple integrals are of the same order and this order is at most 4, we prove that two random variables in the same Wiener chaos either admit a joint density, either are proportional. The basic idea is to write the Malliavin matrix as a sum of squares and to compute the dominant term of its determinant. 

We organized our paper as follows. Section 2 contains some preliminaries on analysis on Wiener chaos. Section 3 is devoted to express the Malliavin matrix  as the sum of the squares of some random variables and in Section 4 we derive an explicit formula for the determinant of $\Lambda$ which also involves the determinant of the covariance matrix. In Section 5 we discuss the existence of the joint density of a vector of multiple integrals.

\section{Preliminaries}

We briefly describe the tools from the analysis on Wiener space that we will need in our work. For complete presentations, we refer to \cite{N} or \cite{NPbook}.
Let $H$ be a real and separable Hilbert space and consider $(W(h), h\in H)$ an isonormal process. That is, $(W(h), h\in H)$ is a family of centered  Gaussian random variables on the probability space $(\Omega, {\cal{F}}, P)$  such that $E W(h) W(g)= \langle f, g\rangle _{H}$ for every $h, g\in H$. Assume that the $\sigma$-algebra ${\cal{F}}$ is generated by $W$. 

Denote,  for $n\geq 0$,  by ${\cal{H}} _{n}$ the $n$th  Wiener chaos generated by $W$. That is, ${\cal{H}} _{n}$ is the vector subspace of $L^{2}(\Omega)$ generated by $\left( H_{n} (W(h)), h\in H, \Vert h\Vert =1\right)$ where $H_{n}$ the Hermite polynomial of degree $n$.  For any $n\geq 1$, the mapping $I_{n}(h^{\otimes n})= H_{n}(W(h))$ can be extended to an isometry between the Hilbert space  $H^{\otimes n}$ endowed with the norm $\sqrt{n!} \Vert \cdot \Vert _{H ^{\otimes n}} $ and the $n$th  Wiener chaos ${\cal{H}}_{n}$. The random variable $I_{n}(f)$ is called the multiple Wiener It\^o integral of $f$ with respect to $W$.

Consider $(e_{j})_{j\geq 1}$ a complete orthonormal system in $H$ and let $f\in H^{\otimes n}$, $ g\in H^{\otimes m}$ be two symmetric functions with $n,m\geq 1$.  Then
\begin{equation}\label{f}
f= \sum_{j_{1}, .., j_{n} \geq 1} \lambda _{j_{1},.., j_{n}} e_{j_{1}} \otimes...\otimes e_{j_{n}}
\end{equation}
and
\begin{equation}\label{g}
g= \sum_{k_{1},.., k_{m}\geq 1} \beta _{ k_{1},.., k_{m}} e_{k_{1}}\otimes..\otimes e_{k_{m}} 
\end{equation}
where the coefficients $\lambda _{i}$ and $\beta _{j}$  satisfy $\lambda_{ j_{\sigma (1)},...j_{\sigma (n)}}= \lambda _{j_{1},.., j_{n}}$ and $ \beta _{ k_{\pi(1)},..., k_{\pi(m)}}= \beta _{ k_{1},.., k_{m}}$ for every permutation $\sigma $ of the set $\{1,..., n\}$ and for every permutation $\pi$ of the set $\{1,.., m\}$. Actually $\lambda _{j_{1},.., j_{n}}  = \langle f, e_{j_{1}} \otimes...\otimes e_{j_{n}}\rangle$ and $\beta _{ k_{1},.., k_{m}} =\langle g, e_{k_{1}}\otimes..\otimes e_{k_{m}} \rangle $ in (\ref{f}) and (\ref{g}).  Note that, throughout the paper we will use the notation $\langle \cdot, \cdot \rangle $ to indicate the scalar product in $H^{\otimes k}$,  independently of  $k$.

If $f\in H^{\otimes n}$, $ g\in H^{\otimes m}$ are symmetric  given by (\ref{f}), (\ref{g}) respectively, then the contraction of order $r$ of $ F$ and $g$ is given by
\begin{eqnarray}
f\otimes _{r}g &=& \sum_{ i_{1},.., i_{r} \geq 1} \sum_{j_{1},..,j_{n-r}\geq 1} \sum_{k_{1},.., k_{m-r}\geq 1} \lambda _{i_{1},.., i_{r} , j_{1}, .., j_{n-r}}\beta _{i_{1},.., i_{r} , k_{1}, .., k_{m-r}}\nonumber\\
&&\times \left( e_{j_{1}}\otimes ..\otimes e_{j_{n-r}}\right) \otimes \left( e_{k_{1}}\otimes ..\otimes e_{k_{m-r}}\right) \label{contra}
\end{eqnarray}
for every $r=0,.., m\wedge n$. In particular $f\otimes _{0}g= f\otimes g$.
Note that $f\otimes _{r} g $ belongs to $H^{\otimes (m+n-2r)}$ for every $r=0,.., m\wedge n
$ and it is not in general symmetric. We will denote by $f\tilde{\otimes }_{r}g$ the symmetrization of $f\otimes _{r}g$.  In the particular case when $H= L^{2}(T, {\cal{B}}, \mu)$ where $\mu$ is a sigma-finite measure without atoms, (\ref{contra}) becomes
\begin{eqnarray}
&&(f\otimes _{r}g) (t_{1},.., t_{m+n-2r})\nonumber\\
&=& \int_{T^{r}}d\mu(u_{1})..d\mu (u_{r}) f(u_{1}, .., u_{r}, t_{1},.., t_{n-r})g(u_{1},..,u_{r}, t_{n-r+1},.., t_{m+n-2r})  \label{contra2}
\end{eqnarray}

An important role will be played by  the following product formula for multiple Wiener-It\^o integrals: if $f\in H^{\otimes n}$, $ g\in H^{\otimes m}$ are symmetric, then
\begin{equation}
\label{prod}
I_{n}(f) I_{m}(g)= \sum_{r=0} ^{m\wedge n} r! C_{m}^{r} C_{n}^{r} I_{m+n-2r} \left(f\tilde{\otimes }_{r}g\right).
\end{equation}
We will need the concept of Malliavin derivative $D$ with respect to $W$, but we will use only its  action  on Wiener chaos. In order to avoid  too many details,  we will just say that, if $f$ is given by (\ref{f}) and $ I_{n}(f)$ denotes its multiple integral of order $n$ with respect to $W$, then 
\begin{equation*}
DI_{n}(f) =n \sum_{j_{1}, .., j_{n} \geq 1} \lambda _{j_{1},.., j_{n}} I_{n-1} \left( e_{j_{2}}\otimes..\otimes e_{j_{n}}\right)e_{j_{1}}.
\end{equation*}
If $F,G$ are two random variables which are differentiable in the Malliavin sense, we will denote throughout the paper by $C$ the covariance matrix and  by  $\Lambda $  the Malliavin matrix of the random vector $(F,G)$. That is, 
\[
\Lambda =\left(  \begin{array}{cc}
\Vert DF\Vert _{H} ^{2}&\langle DF, DG\rangle_{H}\\ \langle DF, DG\rangle_{H} &\Vert DF\Vert _{H} ^{2}
\end{array} \right) .\]

\section{The Malliavin matrix as a sum of squares}

In this section we will express the determinant of the Malliavin matrix of a random couple  as a sum of squares of certain random variables. This will be useful in order to derive the exact formula for the determinant of the Malliavin matrix and its connection with the determinant of the covariance matrix for a given random vector of dimension 2.

Let $f \in H^{\otimes n}$ and $g\in H ^{\otimes m}$ be given by (\ref{f}) and (\ref{g}) respectively, with $n,m\geq 1$.  Let $F=I_{n}(f), G=I_{m}(g) $ denote the multiple Wiener-It\^o integrals of $f$ and $g$   with respect to $W$ respectively. Then
\begin{equation}\label{if}
I_{n} (f)= \sum_{j_{1}, .., j_{n} \geq 1} \lambda _{j_{1},.., j_{n}}I_{n}\left(  e_{j_{1}} \otimes...\otimes e_{j_{n}}\right)
\end{equation}
and
\begin{equation}\label{ig}
I_{m}(g)= \sum_{k_{1},.., k_{m}\geq 1} \beta _{ k_{1},.., k_{m}}I_{m}\left( e_{k_{1}}\otimes..\otimes e_{k_{m}} \right).
\end{equation}
From (\ref{if})  and (\ref{ig}) we have 
\begin{equation*}
DF =n \sum_{j_{1}, .., j_{n} \geq 1} \lambda _{j_{1},.., j_{n}} I_{n-1} \left( e_{j_{2}}\otimes..\otimes e_{j_{n}}\right)e_{j_{1}}
\end{equation*}
and
\begin{equation*}
DG= m \sum_{k_{1},.., k_{m}\geq 1} \beta _{ k_{1},.., k_{m}}I_{m-1}\left( e_{k_{2}}\otimes..\otimes e_{k_{m}} \right)e_{k_{1}}.
\end{equation*}
This implies
\begin{equation*}
\Vert DF\Vert ^{2}_{H} =n^{2} \sum _{i\geq 1} \sum _{j_{2},.., j_{n} \geq 1}\sum_{k_{1},.., k_{n} \geq 1} \lambda _{i, j_{2},.., j_{n}}\lambda _{i, k_{2},.., k_{n}}I_{n-1} \left( e_{j_{2}}\otimes..\otimes e_{j_{n}}\right)I_{n-1} \left( e_{k_{2}}\otimes..\otimes e_{k_{n}}\right)
\end{equation*}
and
\begin{equation*}
\Vert DG\Vert ^{2}_{H} = m^{2}\sum _{l\geq 1} \sum _{l, j_{2},.., j_{n} \geq 1}\sum_{l, k_{2},.., k_{n} \geq 1} 
\beta _{ l, j_{2},.., j_{m}}\beta _{ l, k_{2},.., k_{m}}I_{m-1} \left( e_{j_{2}}\otimes..\otimes e_{j_{m}}\right)I_{m-1}\left( e_{k_{2}}\otimes..\otimes e_{k_{m}} \right)
\end{equation*}
and 
\begin{equation*}
\langle DF, DG \rangle _{H}= nm \sum _{i\geq 1} \sum _{j_{2},.., j_{n} \geq 1}\sum_{k_{1},.., k_{m} \geq 1} \lambda _{i, j_{2},.., j_{n}}\beta _{ i, j_{1},.., j_{m}}I_{n-1} \left( e_{j_{2}}\otimes..\otimes e_{j_{n}}\right)I_{m-1}\left( e_{k_{2}}\otimes..\otimes e_{k_{m}} \right).
\end{equation*}
Let us make the following notation. For every $i\geq 1$, let 

\begin{equation}
\label{sif}
S_{i,f}=n \sum _{i\geq 1} \sum _{j_{2},.., j_{n} \geq 1}\lambda _{i, j_{2},.., j_{n}}I_{n-1} \left( e_{j_{2}}\otimes..\otimes e_{j_{n}}\right)
\end{equation}
and
\begin{equation}
\label{sig}
S_{i,g} =m\sum _{i\geq 1}\sum_{i, k_{2},.., k_{m} \geq 1}\beta _{ i, k_{2},.., k_{m}}I_{m-1}\left( e_{k_{2}}\otimes..\otimes e_{k_{m}} \right).
\end{equation}
We can write
\begin{equation*}
\label{d1}
\Vert DF\Vert _{ H}^{2} = \sum_{i\geq 1} S_{i,f} ^{2}, \hskip0.3cm \Vert DG\Vert _{H} ^{2}= \sum_{ l\geq 1} S_{l,g} ^{2} ,
\langle DF, DG \rangle =\sum_{i\geq 1} S_{i,f} S_{i, g}
\end{equation*}
and
\begin{eqnarray*}
det (\Lambda )&=& \Vert DF\Vert _{H} ^{2}\Vert DG\Vert _{H} ^{2}-\langle DF, DG\rangle_{H}^{2}= \sum_{i,l\geq 1} S_{i, f} ^{2}S_{l, g} ^{2}  - \left( \sum_{i\geq 1} S_{i,f} S_{i, g} \right) ^{2}.
\end{eqnarray*}
A key observation is that

\begin{eqnarray}\label{key} 
\sum_{i,l\geq 1} S_{i, f} ^{2}S_{l, g} ^{2}  - \left( \sum_{i\geq 1} S_{i,f} S_{i, g} \right) ^{2}
&=& \frac{1}{2}\sum_{i,l\geq 1 } \left( S_{i,f} S_{l, g} - S_{l,f} S_{i,g} \right) ^{2}. 
\end{eqnarray}
We obtained
\begin{prop}\label{t1}
The determinant of the Malliavin matrix $\Lambda $ of the vector $(F,G)= (I_{n}(f), I_{m}(g) ) $ can be expressed as 
\begin{equation*}
det \Lambda = \frac{1}{2}\sum_{i,l\geq 1} \left( S_{i,f} S_{l, g} - S_{l,f} S_{i,g} \right) ^{2}
\end{equation*}
where $S_{i,f}, S_{i,g}$ are given by (\ref{sif}) and (\ref{sig}) respectively.
\end{prop}

\begin{corollary}
The determinant of the Malliavin matrix $\Lambda $ of the vector $(F,G)= (I_{n}(f), I_{m}(g) ) $ can be expressed as
\begin{equation*}
det \Lambda = \frac{1}{2}\sum_{i,l\geq 1} \left(  \langle DF, e_{i} \rangle \langle DG, e_{l}\rangle -
 \langle DF, e_{l} \rangle \langle DG, e_{i}\rangle\right) ^{2}
\end{equation*}
{\bf Proof: } This comes from Proposition \ref{t1} and  the relations
\begin{equation*}
S_{i,f}= \langle DF, e_{i}\rangle, \hskip0.3cm S_{i,g}= \langle DG, e_{i} \rangle
\end{equation*}
for every $i\geq 1$. \qed

\end{corollary}

\section{The determinant of the Malliavin matrix on Wiener chaos}
Fix $n,m\geq 1$ and $f, g$ in $H^{\otimes n}, H^{\otimes m}$ respectively defined by (\ref{f}) and (\ref{g}). Consider the random vector $(F,G)= ( I_{n}(f), I_{m}(g))$ and denote by $\Lambda $ its Malliavin matrix and by $C$ its covariance matrix.  

Let us compute $E\det \Lambda$. Denote, for every $i,l\geq 1$ 
\begin{equation}\label{sf}
s_{i,f}= n\sum_{j_{2},.., j_{m} \geq 1} \lambda _{ i, j_{2},.., j_{n}} e_{j_{2}}\otimes..\otimes  e_{j_{n}}
\end{equation}
and
\begin{equation}\label{sg}
s_{l,g} =m \sum_{k_{2}, .., k_{m} \geq 1} \beta _{l. k_{2},.., k_{m}}e_{k_{2}}\otimes..\otimes  e_{k_{m}}.
\end{equation}
Clearly, for every $i,l\geq 1$
\begin{equation}\label{x1}
S_{i,f}= I_{n-1} (s_{i,f}), \hskip0.3cm S_{i,g}= I_{m-1} (s_{i,g} ).
\end{equation}

The following lemma plays a key role in our construction.
\begin{lemma}\label{l5}If $f\in H ^{\otimes n}$ and $ g\in H ^{\otimes m} $ are given by (\ref{f}) and (\ref{g}) respectively and $s_{i,f}, s_{i,g}$ by (\ref{sf}), (\ref{sg}) respectively, then  for every $r=0,.., n-1$
\begin{equation*}
 f\otimes _{r+1} g =\frac{1}{nm} \sum_{i\geq 1} \left( s_{i,f} \otimes _{r} s_{i,g} \right).
\end{equation*}

\end{lemma}
{\bf Proof: } Consider first $r=0$. Clearly, by (\ref{contra})
\begin{eqnarray*}
f\otimes _{1}g &=& \sum_{i\geq 1} \sum_{j_{2},.., j_{n}\geq 1} \sum _{k_{2}, .., k_{m}\geq 1} \lambda _{i,j_{2},.., j_{n}} \beta_{i, k_{2}, .., k_{m}} e_{j_{2}}\otimes..\otimes e_{j_{n}}\otimes e_{k_{2}}\otimes ..e_{k_{m}}\\
&=&\frac{1}{nm} \sum_{i\geq 1} \left( s_{i,f} \otimes s_{i,g} \right).
\end{eqnarray*}
The same argument applies for every $r=1,.., n-1$. Indeed, 
\begin{eqnarray*}
&&f\otimes_{r+1} g\\
&=&\left(  \sum_{j_{1},.., j_{n} \geq 1} \lambda _{j_{1},.., j_{n}} e_{j_{1}}\otimes..\otimes e_{j_{n}}\right)
\otimes_{r}  \left( \sum_{ k_{1},.., k_{m}\geq 1} \beta _{k_{1},.., k_{m}} e_{k_{1}}\otimes..\otimes e_{k_{m}} \right) \\
&=&\sum_{ i_{1},..,i_{r+1}} \sum_{j_{r+2},.., j_{n}} \sum_{k_{r+2},.., k_{m}} \lambda  _{i_{1},..,i_{r+1}, j_{r+2}, .., j_{n}}
\beta _{i_{1},..,i_{r+1}, k_{r+2}, .., k_{m}}\left( e_{j_{r+2}}\otimes..\otimes e_{j_{n}}\right) \otimes \left( e_{k_{r+2}},.., e_{k_{m}}\right)
\end{eqnarray*}
and by (\ref{contra}) again
\begin{eqnarray*}
&&\sum_{i\geq 1} s_{i,j}\otimes _{r} s_{i,g} \\
&=&nm \sum_{i\geq 1} \left( \sum_{j_{2},.., j_{n}\geq 1} \lambda _{i, j_{2},.., j_{n}}e_{j_{2}}\otimes..\otimes e_{j_{n}}\right)\otimes _{r} \left( \sum_{ k_{2},.., k_{m}\geq 1} \beta _{i, k_{2},.., k_{m}} e_{k_{2}}\otimes..\otimes e_{k_{m}}\right)\\
&=&nm \sum_{i\geq 1}\sum_{i_{2},..,i_{r+1}} \sum_{j_{r+2},.., j_{n}} \sum_{k_{r+2},.., k_{m}} 
\lambda  _{i, i_{2},..,i_{r+1}, j_{r+2}, .., j_{n}}
\beta _{i, i_{2},..,i_{r+1}, k_{r+2}, .., k_{m}}\\
&&\times \left( e_{j_{r+2}}\otimes..\otimes e_{j_{n}}\right) \otimes \left( e_{k_{r+2}},.., e_{k_{m}}\right)\\
&=&nm\sum_{ i_{1},..,i_{r+1}} \sum_{j_{r+2},.., j_{n}} \sum_{k_{r+2},.., k_{m}} \lambda  _{i_{1},..,i_{r+1}, j_{r+2}, .., j_{n}}\beta _{i_{1},..,i_{r+1}, k_{r+2}, .., k_{m}}\left( e_{j_{r+2}}\otimes..\otimes e_{j_{n}}\right) \otimes \left( e_{k_{r+2}},.., e_{k_{m}}\right)\\
&=& f\otimes_{r+1} g.
\end{eqnarray*}\qed

\vskip0.2cm

We make a first step to compute $E \det \Lambda$.

\begin{lemma}\label{l7}Let $f\in H ^{\otimes n}, g\in H ^{\otimes m}$ be symmetric and denote by $\Lambda $ the Malliavin matrix of the vector $(F,G)= (I_{n}(f), I_{m}(g))$. Then we have
\begin{equation*}
E \det \Lambda = \sum_{k=0} ^{(n-1) \wedge (m-1)}  T_{k}
\end{equation*}
where we denote, for $k=0,..,(m-1)\wedge (n-1)$, 
\begin{equation}
\label{tk}
T_{k}:=  \frac{1}{2}\sum_{i,l\geq 1}k!^{2} \left( C_{m-1}^{k} \right) ^{2} \left(C_{m-1}^{k} \right) ^{2} (m+n-2-2k)! \Vert s_{i,f}\tilde{\otimes }_{k} s_{l,g} -s_{l,f}\tilde{\otimes }_{k} s_{i,g}\Vert  ^{2}
\end{equation}and $s_{i,f}, s_{i,g}$ are given by (\ref{sf}), (\ref{sg}) for $i\geq 1$. 
\end{lemma}
{\bf Proof: } By Proposition \ref{t1} and relation (\ref{x1})
\begin{eqnarray*}
2 \det \Lambda &=& \sum_{i,l\geq 1} \left( I_{n-1} (s_{i,f}) I_{m-1} (s_{l,g}) - I_{n-1} (s_{l,f}) I_{m-1} (s_{i,g})  \right) ^{2}\\ 
&=&\sum_{i,l\geq 1} \left(\sum_{k=0} ^{(m-1)\wedge (n-1)} k! C_{m-1}^{k} C_{n-1} ^{k}I_{m+n-2-2k}\left( s_{i,f}\tilde{\otimes _{k}} s_{l,g} -s_{l,f}\tilde{\otimes _{k}} s_{i,g}\right)\right) ^{2}
\end{eqnarray*}
where we used the  the product formula  (\ref{prod}). Consequently, from the isometry of multiple stochastic integrals,
\begin{eqnarray*}
E\det \Lambda &=& \frac{1}{2}\sum_{i, l\geq 1} \sum_{k=0} ^{(n-1) \wedge (m-1)} k!^{2} \left( C_{m-1}^{k} \right) ^{2} \left(C_{n-1}^{k} \right) ^{2} (m+n-2-2k)! \Vert s_{i,f}\tilde{\otimes }_{k} s_{l,g} -s_{l,f}\tilde{\otimes }_{k} s_{i,g}\Vert  ^{2}\\
&=& \sum_{k=0} ^{(n-1) \wedge (m-1)}  T_{k}.
\end{eqnarray*}
\qed \vskip0.2cm

For every $n,m\geq1$  let us denote by
\begin{equation}
\label{rnm}
R_{n,m}:=  \sum_{k=1} ^{(n-1) \wedge (m-1)}  T_{k}, \hskip0.6cm R_{n}:=R_{n,n}.
\end{equation}

\begin{remark}
Obviously all the terms $T_{k}$ above are positive, for $k=0,..,(n-1)\wedge (n-1)$. 
\end{remark}

We will need two more auxiliary lemmas.
\begin{lemma}\label{l3}

\ Assume $f_{1}, f_{3} \in H^{\otimes n}$ and $f_{2}, f_{4} \in H^{\otimes m}$ are symmetric functions. Then for every $r=0,.., (m-1)\wedge (n-1)$ we have
\begin{equation*}
\langle f_{1} \otimes _{n-r}f_{3} ,  f_{2} \otimes _{m-r}f_{4} \rangle= \langle f_{1} \otimes_{r} f_{2}, f_{3}\otimes _{r} f_{4} \rangle.
\end{equation*}
\end{lemma}
{\bf Proof: } The case $r=0$ is trivial, so assume $r\geq 1$. Without any loss of the generality, assume that $H$ is $L^{2}(T; \mu) $ where $\mu $ is a sigma-finite measure without atoms. Then, by (\ref{contra2})
\begin{eqnarray*}
&&\langle f_{1} \otimes _{n-r}f_{3} ,  f_{2} \otimes _{m-r}f_{4} \rangle\\
&& \int_{T^{r} }d\mu ^{r}(t_{1},.., t_{r} ) \int_{T^{r}} d\mu ^{r}(s_{1},.., s_{r} )  \\
&&\left( \int_{ T^{n-r}} d\mu^{n-r}(u_{1},..,u_{n-r})f_{1}(u_{1},..,u_{n-r},t_{1},..,t_{r} ) f_{3}(u_{1},.., u_{n-r}, s_{1},..,s_{r} ) \right)\\
&&\left( \int_{ T^{m-r}} d\mu^{m-r}(v_{1},..,v_{m-r})f_{2}(v_{1},..,v_{m-r}, t_{1},.., t_{r})f_{4}(v_{1},..,v_{m-r}, s_{1},..,s_{r})\right)\\
&=& \int_{ T^{n-r}} d\mu^{n-r}(u_{1},..,u_{n-r}) \int_{ T^{m-r}} d\mu^{m-r}(v_{1},..,v_{m-r})\\
&&(f_{1}\otimes_{r}f_{2})(u_{1},.., u_{n-r}, v_{1},.., v_{m-r} ) (f_{3}\otimes _{r} f_{4} )(u_{1},.., u_{n-r}, v_{1},.., v_{m-r} )\\
&=&\langle f_{1} \otimes_{r} f_{2}, f_{3}\otimes _{r} f_{4} \rangle.
\end{eqnarray*}
\qed

\begin{lemma}\label{l1}
Suppose $f_{1}, f_{4}\in H ^{\otimes n}, f_{2}, f_{3} \in H^{\otimes m}$ are symmetric functions. Then 
\begin{equation*}
\langle f_{1} \tilde{\otimes }f_{2}, f_{3} \tilde{\otimes} f_{4} \rangle = \frac{m! n!}{(m+n)! }\sum_{r=0} ^{m\wedge n} C_{n}^{r} C_{m} ^{r} \langle f_{1}\otimes _{r} f_{3} , f_{4} \otimes _{r}f_{2} \rangle .
\end{equation*}

\end{lemma} 
{\bf Proof: } This has been stated and proven in \cite{NoRo} in the case $m=n$. Exactly the same lines of the proofs apply for $m\not=n$. \qed \vskip0.2cm

We first compute the term $T_{0} $ obtained for $k=0$ in (\ref{tk}). 

\begin{prop}
\label{pt0} Let $T_{0} $ be given by  (\ref{tk}) with $k=0$. 
\begin{eqnarray*}
T_{0}&=&\sum_{r=0}^{ (n-1)\wedge (m-1)} mn m!n!C_{n-1}^{r} C_{m-1}^{r}\left[  \Vert f\otimes _{r} g\Vert ^{2} -  \Vert f\otimes _{r+1} g\Vert ^{2}\right].
\end{eqnarray*}
\end{prop}
{\bf Proof: } From (\ref{tk}),
\begin{eqnarray*}
T_{0}&=& \frac{1}{2}(m+n-2)!\sum_{i,l \geq 1} \Vert  s_{i,f} \tilde{ \otimes} s_{l, g} - s_{l,f} \tilde{ \otimes} s_{i, g} \Vert ^{2}\\
&=&\frac{1}{2}
(m+n-2)!\sum_{i,l \geq 1} \left[ \Vert  s_{i,f} \tilde{ \otimes} s_{l, g} \Vert ^{2}+ \Vert  s_{l,f} \tilde{ \otimes} s_{i, g} \Vert ^{2}-2 \langle s_{i,f} \tilde{ \otimes} s_{l, g}, s_{l,f} \tilde{ \otimes} s_{i, g}\rangle \right].
\end{eqnarray*}

Let us apply Lemma \ref{l1}  to compute these norms and scalar products.   We  obtain, by letting $f_{1}= s_{i,f}= f_{4}$ and $f_{2}= s_{l,g}= f_{3}$ (note that $s_{i,f}, s_{i,g}$ are symmetric functions in $H^{\otimes n}, H^{\otimes m}$ respectively)

 \begin{eqnarray*}
(m+n-2)! \langle s_{i,f} \tilde{ \otimes} s_{l, g} , s_{i,f} \tilde{ \otimes} s_{l, g} \rangle 
&=& (m+n-2)! \langle s_{i,f} \tilde{ \otimes} s_{l, g} , s_{l,g} \tilde{ \otimes} s_{i, f} \rangle \\
&=& (m-1)! (n-1)!  \sum_{r=0}^{ (n-1)\wedge (m-1)}C_{n-1}^{r} C_{m-1}^{r} \langle s_{i,f} \otimes _{r} s_{l,g} , s_{i,f} \otimes _{r} s_{l,g} \rangle \\
&=&(m-1)! (n-1)!  \sum_{r=0}^{ (n-1)\wedge (m-1)}C_{n-1}^{r} C_{m-1}^{r} \Vert s_{i,f}\otimes _{r} s_{l,g}\Vert ^{2}. 
\end{eqnarray*}
Analogously, for $f_{1}= s_{l,f}=f_{4}$ and $f_{2}= s_{i,g}=f_{3} $ in Lemma \ref{l1} we get
\begin{eqnarray*}
&&(m+n-2)! \langle s_{l,f} \tilde{ \otimes} s_{i, g} , s_{l,f} \tilde{ \otimes} s_{i, g} \rangle \\
&&= \sum_{r=0}^{ (n-1)\wedge (m-1)} (n-1)! (m-1)!C_{n-1}^{r} C_{m-1}^{r} \Vert s_{l,f} \otimes _{r} s_{i,g} \Vert ^{2}.
\end{eqnarray*}
Next, with $f_{1}= s_{i,f}, f_{2}= s_{l,g}, f_{4}= s_{l,f}, f_{3}= s_{i,g}$
\begin{eqnarray*}
&&(m+n-2)! \langle s_{i,f} \tilde{ \otimes} s_{l, g} , s_{l,f} \tilde{ \otimes} s_{i, g} \rangle \\
&=&(m+n-2)! \langle s_{i,f} \tilde{ \otimes} s_{l, g} , s_{i,g} \tilde{ \otimes} s_{l, f} \rangle \\
&=&\sum_{r=0}^{ (n-1)\wedge (m-1)} (n-1)! (m-1)!C_{n-1}^{r} C_{m-1}^{r} \langle s_{i,f} \otimes _{r} s_{i,g} , s_{l,f} \otimes _{r} s_{l,g}\rangle.
\end{eqnarray*}
Then
\begin{eqnarray}
&&(m+n-2)!\sum_{i,l\geq 1} \Vert s_{i,f} \tilde{ \otimes} s_{l, g} - s_{l,f} \tilde{\otimes } s_{i,g} \Vert ^{2}\nonumber \\
&=&  \sum_{r=0}^{ (n-1)\wedge (m-1)} (n-1)! (m-1)!C_{n-1}^{r} C_{m-1}^{r}\nonumber\\
&&\sum_{i,l\geq 1} \left[ \Vert s_{l,f} \otimes _{r} s_{i,g} \Vert ^{2}+\Vert s_{i,f} \otimes _{r} s_{l,g} \Vert ^{2}
-2 \langle s_{i,f} \otimes _{r} s_{i,g} , s_{l,f} \otimes _{r} s_{l,g}\rangle\right]\nonumber\\
&=&
2 \sum_{r=0}^{ (n-1)\wedge (m-1)} (n-1)! (m-1)!C_{n-1}^{r} C_{m-1}^{r}\sum_{i,l\geq 1}\left[  \Vert s_{i,f} \otimes _{r} s_{l,g} \Vert ^{2}-\langle s_{i,f} \otimes _{r} s_{i,g} , s_{l,f} \otimes _{r} s_{l,g}\rangle \right]\nonumber\\
&=&2 \sum_{r=0}^{ (n-1)\wedge (m-1)} (n-1)! (m-1)!C_{n-1}^{r} C_{m-1}^{r} \nonumber\\
&&\times \left[ \sum_{i,l\geq 1} \Vert s_{i,f} \otimes _{r} s_{l,g} \Vert ^{2}-\langle \sum _{i\geq 1} s_{i,f}\otimes _{r} s_{i,g}, \sum_{l\geq 1} s_{l,f} \otimes _{r} s_{l,g} \rangle\right] \nonumber\\
&=&2 \sum_{r=0}^{ (n-1)\wedge (m-1)} (n-1)! (m-1)!C_{n-1}^{r} C_{m-1}^{r} \left[ \sum_{i,l\geq 1} \Vert s_{i,f} \otimes _{r} s_{l,g} \Vert ^{2}-\Vert \sum _{i\geq 1} s_{i,f}\otimes _{r} s_{i,g}\Vert ^{2}\right].\label{u3}
\end{eqnarray}
Notice that, by Lemma \ref{l5}, for every $r=0,.., n-1$
\begin{equation}
\label{u1}
\Vert \sum _{i\geq 1} s_{i,f}\otimes _{r} s_{i,g}\Vert ^{2}=n^{2}m^{2}\Vert f\otimes _{r+1}g\Vert ^{2}. 
\end{equation}
We apply now Lemma \ref{l3}   and we get 
\begin{eqnarray*}
\sum_{i,l\geq 1} \Vert s_{i,f} \otimes _{r} s_{l,g} \Vert ^{2}
&=&\sum_{i,l\geq 1} \langle  s_{i,f} \otimes _{r} s_{l,g}, s_{i,f} \otimes _{r} s_{l,g}\rangle \\
&=& \sum_{i,l\geq 1} \langle s_{i,f} \otimes_{n-1-r} s_{i,f}, s_{l,g}\otimes_{m-1-r} s_{l,g} \rangle \\
&=&\langle \sum_{i\geq 1}   \langle s_{i,f} \otimes_{n-1-r} s_{i,f}, \sum_{l\geq 1} s_{l,g} \otimes _{m-r-1}s_{l,g}\rangle 
\end{eqnarray*}
and by Lemma \ref{l5} and Lemma \ref{l3}, this equals
\begin{eqnarray}
\sum_{i,l\geq 1} \Vert s_{i,f} \otimes _{r} s_{l,g} \Vert ^{2}
&=& n^{2}m^{2}\langle f \otimes _{n-r} f, g\otimes _{m-r} g\rangle \nonumber \\
&=& n^{2}m^{2}\Vert f\otimes _{r} g\Vert ^{2}.\label{u2}
\end{eqnarray}
By replacing (\ref{u1}) and (\ref{u2}) in (\ref{u3}) we obtain
\begin{eqnarray*}
T_{0}&=&\frac{1}{2}(m+n-2)!\sum_{i,l\geq 1} \Vert s_{i,f} \tilde{ \otimes} s_{l, g} - s_{l,f} \tilde{\otimes } s_{i,g} \Vert ^{2}\nonumber \\
&=&\sum_{r=0}^{ (n-1)\wedge (m-1)} mn m!n!C_{n-1}^{r} C_{m-1}^{r}\left[  \Vert f\otimes _{r} g\Vert ^{2} -  \Vert f\otimes _{r+1} g\Vert ^{2}\right].
\end{eqnarray*}\qed

As  a consequence of the above proof, we obtain 

\begin{corollary}
For every $r=0,.., (m-1)\wedge (n-1)$ and if $s_{i,f}, s_{i,g}$ are given by (\ref{sf}), (\ref{sg}),  it holds that
\begin{equation*}
n^{2} m^{2} \sum_{r=0} ^{n-1} C_{m-1}^{r} C_{m-1}^{r} \left[ \Vert f\otimes _{r} g \Vert ^{2} -\Vert f\otimes _{r+1} g \Vert ^{2}\right]= \sum_{i,l\geq 1} \Vert  s_{i,f}  \tilde{\otimes }s_{l, g} - s_{l,f} \tilde{\otimes }  s_{i,g} \Vert ^{2}.
\end{equation*}
As a consequence, for every $r=0,.., (m-1)\wedge (n-1)$  we have
\begin{equation}
\label{!}
\sum_{r=0} ^{n-1} C_{m-1}^{r} C_{m-1}^{r} \left[ \Vert f\otimes _{r} g \Vert ^{2} -\Vert f\otimes _{r+1} g \Vert ^{2}\right]\geq 0.
\end{equation}
\end{corollary}
{\bf Proof: } It is a consequence of the proof of Proposition \ref{pt0}. \qed \vskip0.2cm

Let us state the main results of this section.

\begin{theorem}\label{t3}
Let $f\in H ^{\otimes n}, g\in H ^{\otimes m}$($n,m\geq 1$) be symmetric and denote by $\Lambda $ the Malliavin matrix of the vector $(F,G)= (I_{n}(f), I_{m}(g))$. Then 
\begin{equation*}
\det \Lambda = \sum_{r=0}^{ (n-1)\wedge (m-1)} mn m!n!C_{n-1}^{r} C_{m-1}^{r}\left[  \Vert f\otimes _{r} g\Vert ^{2} -  \Vert f\otimes _{r+1} g\Vert ^{2}\right]+ R_{n,m}
\end{equation*}
where for every $n,m\geq 1$, $R_{n,m}$ is given by (\ref{rnm}). Note that $R_{n,m}\geq 0$ for every $n,m\geq 1$. 
\end{theorem}
{\bf Proof: } It follows from Proposition \ref{pt0} and Lemma \ref{l7}. \qed 
\vskip0.3cm

In the case when the two multiple integrals live in the same Wiener chaos, we have a nicer expression. 

\begin{theorem}\label{t2}  Under the same assumptions as in Theorem \ref{t3} but with $m=n$, we have
\begin{equation*}
\det \Lambda = m^{2} \det C+ (mm!)^{2}  \sum_{r=1} ^{ \left[ \frac{m-1}{2}\right]} \left( (C_{m-1}^{r}) ^{2}- (C_{m-1} ^{r-1} ) ^{2} \right) \left( \Vert f\otimes _{r}g \Vert ^{2} -\Vert f\otimes _{n-r} g\Vert ^{2} \right) + R_{m}
\end{equation*}
with $R_{m}$ given by (\ref{rnm}).  Here $[x]$ denotes the integer part of $x$.

\end{theorem} 
{\bf Proof: } Suppose $n\leq m$ and that $m$ is odd.  The case $m$ even is similar. From Theorem \ref{t3} we have
\begin{eqnarray*}
\det \Lambda &=& (m m! ) ^{2}\left[ \sum_{r=0} ^{(m-1)} \left( C_{m-1}^{r}\right) ^{2} \Vert f\otimes _{r} g\Vert ^{2} -\sum_{r=0} ^{(m-1)} \left( C_{m-1}^{r}\right) ^{2} \Vert f\otimes _{r+1} g\Vert ^{2} \right]\\
&=&(m m! ) ^{2}\left[\sum_{r=0} ^{\frac{m-1}{2}} \left( C_{m-1}^{r}\right) ^{2} \Vert f\otimes _{r} g\Vert ^{2} -\sum_{r=\frac{m-1}{2}} ^{(m-1)} \left( C_{m-1}^{r}\right) ^{2} \Vert f\otimes _{r+1} g\Vert ^{2}\right] \\
&&+ (m m! ) ^{2}\left[\sum_{\frac{m-1}{2}}^{m-1} \left( C_{m-1}^{r}\right) ^{2} \Vert f\otimes _{r} g\Vert ^{2} -\sum_{r=0}^{\frac{m-1}{2}} \left( C_{m-1}^{r}\right) ^{2} \Vert f\otimes _{r+1} g\Vert ^{2}\right] \\
&=&  (m m! ) ^{2}\sum_{r=0} ^{\frac{m-1}{2}} \left( C_{m-1}^{r}\right) ^{2}\left[  \Vert f\otimes _{r} g\Vert ^{2}- \Vert f\otimes_{n-r} g\right]\\
&&+(m m! ) ^{2}\sum_{r=1}^{\frac{m-1}{2}} \left( C_{m-1}^{r}\right) ^{2}\left[  \Vert f\otimes _{n-r} g\Vert ^{2}- \Vert f\otimes_{r} g\Vert ^{2}\right]
\end{eqnarray*}
where we made the change of index $r'= n-1-r$ in the second and third sum above. Finally, noticing that for $r=0$ we have 
$$m ^{2} m! ^{2} \left( C_{m-1}^{0}\right) ^{2}\left[  \Vert f\otimes _{0} g\Vert ^{2}- \Vert f\otimes_{n} g\Vert \right]= m ^{2}\det C$$
we obtain the conclusion.  
  \qed

\begin{example}\label{exx}
Suppose  $m=n=2$. Then 
\begin{eqnarray*}
\det \Lambda &=& 16 \left[ \Vert f\otimes g\Vert ^{2} - \Vert f \otimes _{2} g \Vert ^{2} \right] + R_{2} \\
&=& 4\det C + R_{2} .
\end{eqnarray*}
We retrieve the formula in \cite{NoNuPo} with $R_{2} = 32 \left( \Vert f\otimes _{1}g \Vert ^{2}- \Vert f\tilde{\otimes }_{1}g\Vert ^{2} \right). $

Assume $m=n=3$.  Then 
\begin{eqnarray*}
\det \Lambda &=& 9\times 36  \left[ \left( \Vert f\otimes g\Vert ^{2}- \Vert f\otimes _{3} g\Vert ^{2} \right) + 9 \times 36 \times \left( (C_{2}^{1})^{2} -1\right) \left( \Vert f\otimes _{1}g\Vert ^{2}- \Vert f\otimes _{2}g\Vert ^{2}\right) \right] +R_{3}\\
&=& 9\det C + 9 \times 36 \times 3  \left( \Vert f\otimes _{1}g\Vert ^{2}- \Vert f\otimes _{2}g\Vert ^{2}\right)  + R_{3}.
\end{eqnarray*}

Suppose $m=n=4$. Then 
\begin{eqnarray*}
\det \Lambda &=&16 \det C + 16 \times 4!  \times 4!  \left( (C_{3}^{1})^{2} -1\right) \left( \Vert f\otimes _{1}g\Vert ^{2}- \Vert f\otimes _{3}g\Vert ^{2}\right) + R_{4}.
\end{eqnarray*}
\end{example} 

\section{Densities of vectors of multiple integrals}\label{density}

Let us discuss when a couple of multiple stochastic integrals has a law which is absolutely continuous with respect to the Lebesque measure. The situations when the components of the vector are in the same chaos of in chaoses of different orders need to be separated.

Let us first discuss the case of variables in the same chaos. In order better understand the relation between $\det \Lambda $ and $\det C$ we need more information on the terms $R_{m}$ in Theorem \ref{t2}. It is actually possible  to compute the last term $T_{m-1}$ in (\ref{tk}).  

\begin{prop}
Suppose $m=n$ and let $T_{m-1}$ be the term obtained in (\ref{tk}) for $r=m-1$. Then 
\begin{equation*}
T_{m-1}=m^{2} m! ^{2} \left[ \Vert f \otimes _{m-1} g\Vert ^{ 2} -  \langle f\otimes _{1} g, g\otimes_{ 1} f\rangle \right].
\end{equation*}
\end{prop}
{\bf Proof: } From (\ref{tk}),
\begin{eqnarray*}
T_{m-1}&=& \frac{1}{2} \sum_{i, l \geq 1} (m-1) ! ^{2} \Vert | s_{i,f} \tilde{\otimes}_{ m-1} s_{l,g} -s_{l,f} \tilde{\otimes}_{ m-1} s_{i,g}\Vert ^{2} \\
&=& \frac{1}{2} \sum_{i, l \geq 1} (m-1) ! ^{2} \Vert | s_{i,f} \otimes_{ m-1} s_{l,g} -s_{l,f} \otimes_{ m-1} s_{i,g}\Vert ^{2}\\
&=&  \frac{1}{2} \sum_{i, l \geq 1} (m-1) ! ^{2} \left[ \langle s_{i,f}, s_{l,g} \rangle - \langle s_{l,f}, s_{i,g} \rangle\right] ^{2} \\
&=&  (m-1) ! ^{2}\left[ \sum_{i,l\geq 1}  \langle s_{i,f}, s_{l,g} \rangle ^{2} - \sum_{i,l\geq 1}
\langle s_{i,f}, s_{l,g} \rangle\langle s_{l,f}, s_{i,g} \rangle \right]\\
&=& (m-1)! ^{2} \left[ \sum_{i,l \geq 1} \langle s_{i,f} \otimes s_{i,f}, s_{l,g} \otimes s_{l,g} \rangle -  \sum_{i,l \geq 1} \langle s_{i,f} \otimes s_{i,g}, s_{l,g} \otimes s_{l,f} \rangle\right]\\
&=&(m-1)! ^{2}  m^{4} \left[  \langle f\otimes _{1}f, g\otimes _{1} g \rangle - \langle f\otimes _{1} g, g\otimes
_{ 1} f\rangle \right]\\
&=&
m^{2} m! ^{2} \left[ \Vert f \otimes _{m-1} g\Vert ^{ 2} -  \langle f\otimes _{1} g, g\otimes_{ 1} f\rangle \right]
\end{eqnarray*}
where we applied Lemmas \ref{l3} and \ref{l5}.
\qed.

We first answer the open problem 6.2 in \cite{NoNuPo} for chaoses of order lesser than five. 

\begin{theorem}
Let $m\leq 4$ and let $f,g \in H^{\otimes m}$ be symmetric. Then the random vector $(F,G)= (I_{m}(f), I_{m} (g))$ does not admit a density if and only if 
$$\det C =0.$$
 In other words, the vector $(F, G)$ does not admit a density if and only if its components are proportional.

\end{theorem}
{\bf Proof: } The case $m=n=1$ is obvious and the case $m=n=2$ follows from \cite{NoNuPo} (it also follows from Example  \ref{exx}). Suppose $m=n=3$. Then

\begin{eqnarray*}
\det \Lambda &=& 9\det C+ 9\times 36 \times \left( (C_{2}^{1}) ^{2}-1\right) \left[ \Vert f\otimes _{1} g \Vert   ^{2}- \Vert f \otimes _{2} g\Vert  ^{2} \right]  \\
&&+ 9 \times 36 \left[ \Vert f \otimes _{2} g\Vert ^{2} - \langle f\otimes _{2} g, g \otimes _{2} f\rangle \right]  + R_{3}'
\end{eqnarray*}
where $R_{3}'$ is the term with $k=1$ in (\ref{tk}).  Using $\langle f\otimes _{1} g, g \otimes _{1} f\rangle =\langle f\otimes _{2} g, g \otimes _{2} f\rangle$ (Lemma \ref{l3})  we get 
\begin{eqnarray*}
\det \Lambda &=& 9\det C+ 9\times 36 \times 3 \left[ \Vert f\otimes _{1} g \Vert   ^{2} -\langle f\otimes _{1} g, g \otimes _{1} f\rangle \right]  \\
&&- 9\times 36 \times 2 \left[ \Vert f\otimes _{2} g \Vert   ^{2} -\langle f\otimes _{2} g, g \otimes _{2} f\rangle \right]+ R_{3} ' .
\end{eqnarray*}
Suppose  $\det \Lambda =0$. Then $T_{0}, T_{1}, T_{2}$ from (\ref{tk}) vanish. In particular  $T_{2}=0$ in (\ref{tk}) and so  
$$\Vert f\otimes _{2} g \Vert   ^{2} -\langle f\otimes _{2} g, g \otimes _{2} f\rangle=0.$$
This implies 
$$9\det C+ 9\times 36 \times 3 \left[ \Vert f\otimes _{1} g \Vert   ^{2}- \langle f\otimes _{1} g, g \otimes _{1} f\rangle \right] =0$$ 
and therefore $\det C=0$ because $\Vert f\otimes _{1} g \Vert   ^{2} -\langle f\otimes _{1} g, g \otimes _{1} f\rangle $ is positive by Cauchy-Schwarz.

Suppose $m=n=4$. 

\begin{eqnarray*}
\det \Lambda &=&16 \det C+ 16 \times 4! ^{2}  \left( (C_{3}^{1} ) ^{2} -1 \right) \left[ \Vert f\otimes _{1} g \Vert ^{2} - 
\Vert f\otimes _{3} g \Vert ^{2}\right]\\
&& + 16 \times 4! ^{2} \left[ \Vert f\otimes _{3} g\Vert ^{2} - \langle f \otimes _{3} g, g\otimes _{3} f\rangle \right] + R_{4} '
\end {eqnarray*}
where $ R_{4}'$ is the sum of terms obtained for $k=1$ and $ k=2$ in (\ref{tk}).  Since $\langle f \otimes _{3} g, g\otimes _{3} f\rangle=\langle f \otimes _{1} g, g\otimes _{1} f\rangle$(Lemma \ref{l3}) we get
\begin{eqnarray*}
\det \Lambda &=&16 \det C+ 16 \times 4! ^{2}  \left( (C_{3}^{1} ) ^{2} -1 \right) \left[ \Vert f\otimes _{1} g \Vert ^{2} -\langle f \otimes _{1} g, g\otimes _{1} f\rangle\right] \\
&&- 16 \times 4!^{2} \left( (C_{3}^{1} ) ^{2} -2 \right)\left[ \Vert f\otimes _{3} g \Vert ^{2} -\langle f \otimes _{3} g, g\otimes _{3} f\rangle\right]+ R_{4}'.
\end{eqnarray*}
Assume  $\det \Lambda =0$. Then in particular $T_{3}$ from {\ref{tk}) vanishes. So 
$$\Vert f\otimes _{3} g \Vert ^{2} -\langle f \otimes _{3} g, g\otimes _{3} f\rangle=0$$
and this implies $\det C=0$.

\begin{remark}
For $m=n\geq 5$, we have 
\begin{eqnarray*} 
\det \Lambda &=& 25 \det C + 25 \times 5! ^{2}\left( (C_{4}^{1}) ^{2}-1\right)  \left[ \Vert f\otimes _{1} g \Vert ^{2} -  \Vert f\otimes _{4} g \Vert ^{2}\right] \\
&&+ 25 \times 5! ^{2}\left( (C_{4}^{2}) ^{2}-1\right)  \left[ \Vert f\otimes _{2} g \Vert ^{2} -  \Vert f\otimes _{3} g \Vert ^{2}\right] \\
&&+ 25 \times 5! ^{2} \left[ \Vert f\otimes _{4} g \Vert ^{2}-\langle f \otimes _{4} g, g\otimes _{4} f\rangle\right]+ R_{5}'
\end{eqnarray*} 
If $\det \Lambda =0$ then, since $T_{4}$ vanishes, we get that $\Vert f\otimes _{4} g \Vert ^{2}-\langle f \otimes _{4} g, g\otimes _{4} f\rangle$ vanishes. But this is not enough. We need some additional information  in order to handle the difference $\Vert f\otimes _{2} g \Vert ^{2} -  \Vert f\otimes _{3} g \Vert ^{2}$. One possibility is to look to the terms $ T_{1}, T_{2}, T_{3}$ in (\ref{tk}) but these terms cannot be written in a closed form, since they involve more complicated contractions (some "contractions of contractions"). 
\end{remark}

Let us finish by some comments concerning the case of variables in chaoses of different orders. Consider $(F,G)=(I_{n}(f), I_{m} (g))$ with $n\not= m$. 
 First, let us note that $E\det \Lambda =0$ does not imply $\det C=0$.  This can be viewed by considering the following example. 

\begin{example}\label{ex1}
Take $F= I_{2} (f) $ and $G= I_{2}(h^{\otimes 2})$ where $\Vert h\Vert =1$. In this case
$$ \det C=2 \mbox{ and } \det \Lambda = 0.$$
One can also choose $F= I_{n} (h^{\otimes n})$ and $G=I_{m}(h^{\otimes m})$ with $m\not=n$ and $\Vert h\Vert =1$.
\end{example}

In the case  $(I_{n}(f), I_{1}(g)$)  there is only one term in (\ref{tk}) obtained for $k=0$. It reads
\begin{equation*} 
T_{0}= nn! \left[ \Vert f\otimes _{2}  g\Vert ^{2}- \Vert f\otimes _{1}g \Vert ^{2} \right] .
\end{equation*} and therefore the condition for the existence of the joint density is $\Vert f\otimes _{2}  g\Vert ^{2}- \Vert f\otimes _{1}g \Vert ^{2} >0$. 

The case $(I_{n}(f), I_{2}(g)$) is more complicated and needs new ideas in order to obtain the if and only if condition for the existence of the density of the vector. Even the "last term"in (\ref{tk}) (that is, the term obtained for $k=(m-1)\wedge (n-1)$ cannot be written is a nice form.


\begin{thebibliography}{99}

\bibitem{NoNuPo}
{I. Nourdin, D. Nualart and G. Poly (2012): }{\em Absolute continuity and convergence of densities for random vectors on Wiener chaos. } Electronic Journal of Probability,  18, paper 22, 1-19. Arxiv version available at http://arxiv.org/abs/1207.5115v1. 

\bibitem{NPbook}
{I. Nourdin and G. Peccati (2012): } Normal Approximations with Malliavin Calculus
From Stein's Method to Universality. Cambridge University Press.

\bibitem{NoRo}
{I. Nourdin and I. Rosinski (2012): }{\em Asymptotic independence and limit law for multiple Wiener-It\^o integrals. } To appear in The Annals of Probability.

\bibitem{N}
{D. Nualart (2006): }{\em Malliavin Calculus and Related Topics. Second Edition.  }{Springer New York. }








\end{thebibliography}
\end{document}